\documentclass[review]{elsarticle}
\RequirePackage[top=1.5cm,bottom=1.5cm,left=2.cm,right=2.cm,includehead,includefoot]{geometry}

\usepackage{amsmath,amsfonts,amsthm}

\usepackage{booktabs}
\usepackage{longtable}
\usepackage{setspace}
\usepackage{amssymb}
\usepackage{fancyhdr}
\usepackage{mathrsfs}
\usepackage{amsfonts}
\usepackage{cases}
\usepackage{verbatim}
\usepackage{graphicx}      
\usepackage{epstopdf}
\usepackage{grffile}
\graphicspath{{figures/}}
\usepackage{subfig}
\usepackage{caption}
\usepackage{tikz}

\theoremstyle{remark}

\usepackage[frak=mma]{mathalfa}

\begin{document}

\begin{frontmatter}



\title{An Operator Learning Approach via Function-valued Reproducing Kernel Hilbert Space for Differential Equations}
\tnotetext[mytitlenote]{This work was supported by the National Key R$\&$D Program of China (2020YFA0709800),  the National Natural Science Foundation of China (No. 11901577, 11971481, 12071481), the Natural Science Foundation of Hunan (2021JJ20053, 2020JJ5652), the Science and Technology Innovation Program of Hunan Province (No. 2021RC3082), and the Defense Science Foundation of China (2021-JCJQ-JJ-0538).}


\begin{abstract}
Much recent work has addressed the  solution of a family of partial differential equations (PDEs) by computing the inverse operator map between function spaces. Toward this end, we incorporate function-valued reproducing kernel Hilbert spaces (function-valued RKHS) in our operator learning model. Motivated by recently successful neural operator: Deep operator networks (DeepONets), we use neural networks to parameterize the Hilbert-Schmidt integral operator and propose an architecture based on the representer theorem in function-valued RKHS. Experiments including the advection, KdV, burgers’, and poisson equations  show that the proposed architecture has  better accuracy on nonlinear PDEs and linear  PDEs with a small amount of data than DeepONets. We also show that by learning the mappings between function spaces, the proposed method can find the solution of a high-resolution input after learning from lower-resolution data.
\end{abstract}
\author[nudt]{Kaijun Bao}
\author[nudt]{Xu Qian\corref{ccr}}
\cortext[ccr]{Corresponding author.}
\ead{qianxu@nudt.edu.cn}
\author[nudt]{Ziyuan Liu}
\author[nudt,hpc]{Songhe Song}

\address[nudt]{Department of Mathematics, College of Liberal Arts and Science, National University of Defense Technology, Changsha 410073, P.R.China}
\address[hpc]{State Key Laboratory of High Performance Computing, National University of Defense Technology, Changsha 410073, China}

\end{frontmatter}


\section{Introduction}\label{Introduction}
Partial differential equations (PDEs) are powerful tools for modeling the real world in fields such as aerospace and material and biomolecular dynamics, with great success from microscopic  (quantum, molecular dynamics) to cryoscopic issues (ship engineering). However, two   challenges remain: identifying and formulating the   PDEs appropriate to modeling   a specific problem, and computational efficiency when solving complicated  PDE systems.

Modeling a specific problem and determining the approximate underlying PDE usually requires much prior knowledge,   combined with general conservation laws, such as the law of conservation of energy, to design a predictive model. But, to acquire prior knowledge for complex systems  is usually   expensive, or the system is too simple and lacks information, and solving large-scale nonlinear PDEs requires significant computing resources, which can render simulations infeasible. This is reflected in the traditional numerical method for solving PDEs. Since   the analytical solution of PDEs is difficult, in practical engineering applications, we seek the numerical value of a solution at some discrete points in a region. Classical methods such as the finite difference, finite volume, spectrum, and finite element methods will eventually turn a PDE solution into the solution of linear equations. When the problem becomes complicated,  the numerical method must be accurate enough and the discrete points  dense enough, such as in turbulence problems. The distance between discrete points must usually reach the order of micrometers, but this  requires significant computing resources, which renders the problem unsolvable. If the dimension of the problem increases, we only discretize two points in each dimension (this is already quite sparse), and the total number of discrete points increases exponentially,   a problem also faced by traditional numerical methods, i.e., the curse of dimensionality.

Much recent work has used deep neural networks to solve PDEs, showing efficiency compared with  numerical methods. One issue is to parameterize the solution of a PDE with neural networks \cite{4,5,6,7,8,9,24,25}. Based on   automatic differentiation, we can   calculate differential terms in a PDE. Substituting them in the PDE to construct a loss function, the neural network can be well trained. It is worth noting that this approach is   data-independent, because the loss function is   constructed based on knowledge of the underlying PDE structure, whose representative work   is PINN \cite{4}. The deep rize method \cite{5}  of Weinan E  constructs a loss function via the variational form of the PDE, showing that  a  neural network  can  solve high-dimensional problems through Monte Carlo simulation. However,   the neural network must be retrained when the coefficients associated with the PDEs are changed. Numerical methods also have this problem. Other   work focuses on solving a series of PDEs. One issue is to discretize the coefficients and solutions of PDEs at a specific resolution \cite{15,16,17,18}. Based on the data,  a PDE solution is transformed to the learning of a map between finite Euclidean spaces. Such an approach clearly depends on the discretization size and geometry of the training data, and it is impossible to query solutions at new points. Another approach is  the neural operator \cite{2,3,10,11,12,13,14,20,21,22,23}, which is   closest to the problem we investigate. DeepONets \cite{2} take inputs consisting of coefficients $f$ and locations $x$, and generate a solution $u$ at $x$. Mod-net \cite{20} and MWT \cite{13} use the integral operator $T[f](x)=\int_{D} K(x,y)f(y)dy$ to approximate   solutions and use a neural network to model $K$. We can see that the neural operator method considers the locations $x$, and can solve a series of PDEs, but  is independent of on the resolution.

We treat the problem of learning a PDE as an operator learning problem. For a well-posed PDE, given a coefficient, boundary condition, and source term, there exists a unique solution, which means the PDE can be completely represented by the operator mapping coefficients, boundary conditions, or source terms to solutions. This indicates that solving a PDE can be accomplished by learning from data, but we do not need to know its underlying knowledge. Function-valued reproducing kernel Hilbert spaces \cite{1}, which are induced by a unique operator-valued kernel, offer a special form given by the representer theorem to approximate an operator, and the error between the approximate solution given by this approach and the exact solution can be estimated. This indicates that with an appropriate operator-valued kernel and growth of the data, the approximate solution will converge to the exact one at a fixed rate. Inspired by DeepONets, we take the operator-valued kernel as Hilbert-Schmidt integral operator \cite{1}. By using neural networks to parameterize it and combining with the representer theorem, we introduce the function-valued RKHS-based model, a novel operator learning method. Our main contributions are as follows: i) based on the architecture of DeepONets and the representer theorem, we develop a neural network model that efficiently learns the operator map; ii) we demonstrate the applicability of our model on the one-dimensional dataset of nonlinear burgers' and KdV equations and our model performs better than DeepONets; iii) in linear case with a small amount of data, we show that the proposed model can efficiently learn the operator but DeepONets can not; and iv) we show that our model is a mesh-independent method that can find the solution to a high-resolution input after learning from lower-resolution data

\section{Methodology}\label{sec2}
In this section, we will first transform the solution of PDEs to the learning of the operator between two infinite-dimensional spaces. Then, we propose a neural network architecture based on the representer theorem in function-valued RKHS and Hilbert-Schmidt integral operator.
\subsection{Problem Setting}\label{sec2.1}
We aim to solve PDEs by approximating the operator between two infinite-dimensional spaces. The PDE we consider takes the  form
\begin{equation}\label{a1}
\begin{aligned}
(\mathcal{L}_{\alpha}u)(x)&=f(x), &x &\in D \\
u(x)&=g(x), &x &\in \partial D,
\end{aligned}
\end{equation}
where $\mathcal{L}$ is a differential operator, and $\alpha$ involves terms that determine the governing equation of the PDE. We assume that PDE \eqref{a1} is well-posed, which means that there exists an operator mapping from the space constructed by functions $\alpha,f,g$ to the space constructed by solution $u$. We actually fix two terms in functions $\alpha,f,g$, and the PDE can be considered as the operator mapping the remaining term to the solution. We propose to solve the PDE by constructing a parametric map to approximate this operator. Without loss of generality, we fix functions $\alpha,g$.

Let $\mathcal{G}^\dagger$ be the operator mapping source term $f$ to solution $u$ in PDE \eqref{a1}, and $\mathcal{G}$   the parametric map  taking $f$ and  $\theta \in \Theta$ as input, where $\Theta$ represents some finite-dimensional parameter spaces. We expect there is a $\theta^\dagger \in \Theta$ such that $\mathcal{G}(\cdot,\theta^\dagger)\approx\mathcal{G}^\dagger$.

Let $\mathcal{F}$ and $\mathcal{U}$ be Banach spaces of functions and suppose we have observations $\{f_i,u_i\}^n_{i=1}$, where $f_i \backsim \mu$ is an i.i.d. sequence from the probability measure $\mu$ supported on $\mathcal{F}$, and $u_i$ is the corresponding solution, which belongs to $\mathcal{U}$. To determine   $\theta^\dagger$, a natural framework is to define a cost functional $C:\mathcal{U} \times \mathcal{U} \to \mathbb{R}$ and seek a minimizer of the problem,
\begin{equation}\label{a2}
\mathop{\min}_{\theta \in \Theta}\mathbb{E}_{f\backsim \mu}[C(\mathcal{G}(f,\theta),u)].
\end{equation}

To numerically work with $f,u$, we assume pointwise evaluations of functions. Let $P_k$ be a $k$-point discretization of domain $D$, and assume  observations $f_i|_{P_k},u_i|_{P_k}$ for a finite collection of input-output pairs indexed by $i$.

\subsection{Function-valued RKHS-based Model}\label{sec2.2}
Function-valued RKHS, which has the property that the representer theorem provides a special form for approximating an operator, demonstrates that the generalization error only depends on the reproducing kernel. Then, inspired by the architecture of DeepONets, we take Hilbert-Schmidt integral operator as the reproducing kernel and propose our model.  
\begin{figure}[htbp]
\centering
\begin{tikzpicture}
\node[draw,rectangle,rounded corners,](An) at (2,0) {$A_n$};
\node at (2,1) {$\cdots$};
\node[draw,rectangle,rounded corners,](A2) at (2,2) {$A_2$};
\node[draw,rectangle,rounded corners,](A1) at (2,3) {$A_1$};
\node[draw,fill=blue!30,rectangle,rounded corners,](K2) at (2,4) {$K_2$};
\node[draw,rectangle,rounded corners,](fn) at (4,5) {$f_n$};
\node at (4,6) {$\cdots$};
\node[draw,rectangle,rounded corners,](f2) at (4,7) {$f_2$};
\node[draw,rectangle,rounded corners,](f1) at (4,8) {$f_1$};
\node[draw,rectangle,rounded corners,](tn) at (8,0) {$t_n$};
\node at (8,1) {$\cdots$};
\node[draw,rectangle,rounded corners,](t2) at (8,2) {$t_2$};
\node[draw,rectangle,rounded corners,](t1) at (8,3) {$t_1$};
\node[draw,rectangle,rounded corners,](bn) at (8,5) {$l_n$};
\node at (8,6) {$\cdots$};
\node[draw,rectangle,rounded corners,](b2) at (8,7) {$l_2$};
\node[draw,rectangle,rounded corners,](b1) at (8,8) {$l_1$};
\node[draw,fill=gray!30,rectangle,rounded corners,](ji) at (9.5,4) {$\odot$};
\node[draw,rectangle,rounded corners,](v) at (11.5,4) {$v(x)$};
\node[draw,fill=gray!30,rectangle,rounded corners,](o1) at (4,3) {$\otimes$};
\node[draw,fill=gray!30,rectangle,rounded corners,](o2) at (4,2) {$\otimes$};
\node at (4,1) {$\cdots$};
\node[draw,fill=gray!30,rectangle,rounded corners,](on) at (4,0) {$\otimes$};
\node(y) at (0,1.5) {$y$};
\node(x) at (0,4) {$x$};
\node(f) at (0,6.5) {$f$};
\draw[->] (x)--(K2);
\draw[->] (y)--(K2);
\draw[->] (y)--(An);
\draw[->] (y)--(A2);
\draw[->] (y)--(A1);
\draw[->] (f)--(fn);
\draw[->] (f)--(f2);
\draw[->] (f)--(f1);
\draw[->] (A1)--(o1);
\draw[->] (A2)--(o2);
\draw[->] (An)--(on);
\draw[->] (K2)--(o1);
\draw[->] (K2)--(o2);
\draw[->] (K2)--(on);
\draw[->] (f1)--(b1);
\draw[->] (f2)--(b2);
\draw[->] (fn)--(bn);
\draw[->] (on)--(tn);
\draw[->] (o2)--(t2);
\draw[->] (o1)--(t1);
\draw[dashed]
(1.5,-0.5) rectangle (2.5,3.5);
\draw[dashed]
(3.5,4.5) rectangle (4.5,8.5);
\draw[dashed]
(5.5,-0.5) rectangle (6.5,3.5);
\node[draw,fill=gray!30,rectangle,rounded corners,minimum height=3] at (6,1.5) {$\int$};
\draw[dashed]
(5.5,4.5) rectangle (6.5,8.5);
\node[draw,fill=blue!30,rectangle,rounded corners,minimum height=3] at (6,6.5) {$K_1$};
\draw[dashed]
(7.5,-0.5) rectangle (8.5,3.5);
\draw[dashed]
(7.5,4.5) rectangle (8.5,8.5);
\draw[->] (8.5,1.5)--(ji);
\draw[->] (8.5,6.5)--(ji);
\draw[->] (ji)--(v);
\end{tikzpicture}
\caption{Architecture of function-valued RKHS-based model.}
\label{fig:d}
\end{figure}
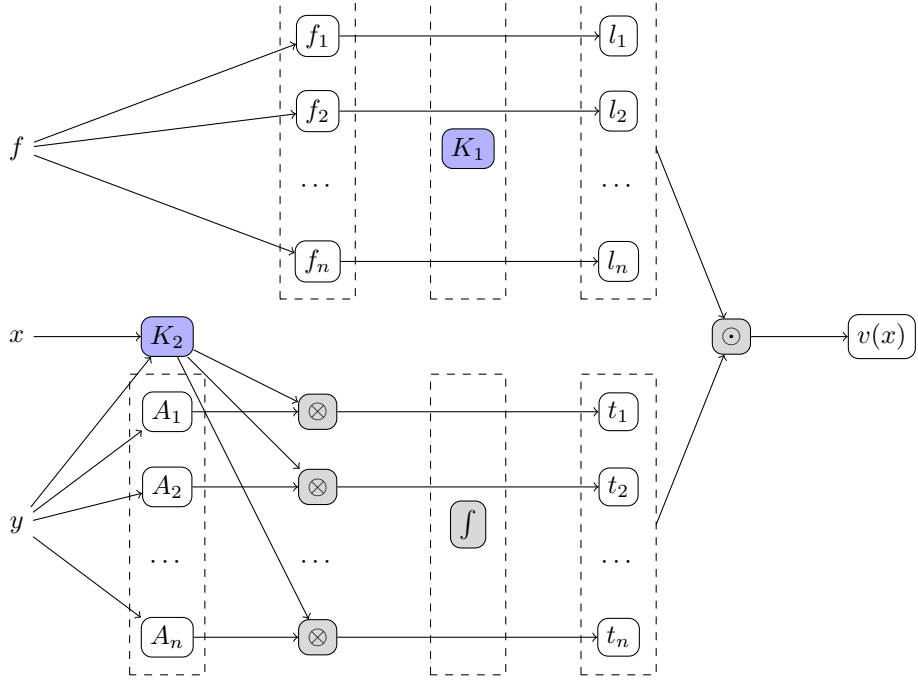

We first consider the supervised learning in view of an operator. Given data $\{f_i,\mathcal{G}^\dagger(f_i)\}^n_{i=1}$, we consider the following regression estimator:
\begin{equation}
\begin{aligned}
\hat{G}_n=\mathop{\rm arg min}\limits_{G \in \mathcal{H}_k}(\hat{\mathcal{R}}_n(G)+\lambda \|G\|^2_{\mathcal{H}_k})  \\
\hat{\mathcal{R}}_n(G)=\frac{1}{n}\sum_{i=1}^n \| G(f_i)-\mathcal{G}^\dagger(f_i) \|^2_{\mathcal{U}},
\end{aligned}
\end{equation}
where $\mathcal{H}_k$ is a Hilbert space. Based on the data, $\hat{G}_n$ is the optimal approximation of $\mathcal{G}^\dagger$, but we wish to know if $\hat{G}_n$ can perform well on the whole input space. The performance is characterized by the generalization error, which takes   the form
\begin{equation}
\mathbb{E}_{f \sim \mu}[\|\hat{G}_n(f)-\mathcal{G}^\dagger(f)\|^2_{\mathcal{U}}].
\end{equation}

Let $\hat{G}=\mathop{\rm arg min}\limits_{G \in \mathcal{H}_k} ({\mathcal{R}}(G)+\lambda \|G\|^2_{\mathcal{H}_k})$, where ${\mathcal{R}}(G)=\mathbb{E}_{f\sim\mu} [\| G(f)-\mathcal{G}^\dagger(f) \|^2_{\mathcal{U}}]$. The generalization error can be decomposed as
\begin{equation}\label{a3}
\mathbb{E}_{f\sim\mu}[\|\hat{G}_n(f)-\mathcal{G}^\dagger(f)\|^2_{\mathcal{U}}] \leqslant \mathbb{E}_{f\sim\mu}[\|\hat{G}_n(f)-\hat{G}(f)\|^2_{\mathcal{U}}]+\mathbb{E}_{f\sim\mu}[\|\hat{G}(f)-\mathcal{G}^\dagger(f)\|^2_{\mathcal{U}}].
\end{equation}

By taking $\mathcal{H}_k$ as the function-valued RKHS, which is induced by a reproducing kernel $K$ and following the representer theorem, we have
\begin{equation}\label{a4}
\hat{G}_n=\sum_{i=1}^n K(f_i,\cdot)A_i,
\end{equation}
where $A_i(x)$ are in $\mathcal{U}$, and   $K$ is a nonnegative operator-valued function. For the first term in equation \eqref{a3}, the approximate theorem in function-valued RKHS shows that with the growth of data, $\hat{G}_n(f)$ can converge to $\hat{G}(f)$ at a fixed rate. Then the second term is what   matters, but it will vanish if $\mathcal{G}^\dagger$ belongs to the function-valued RKHS $\mathcal{H}_k$, which means that it is determined by the reproducing kernel $K$. Hence, based on   equation \eqref{a4}, we parameterize the reproducing kernel $K$ and propose the following neural network architecture.

Inspired by the successful application of DeepONets, we consider the Hilbert-Schmidt integral operator, which has a similar architecture. It takes form:
\begin{equation}
  K(f_i,\cdot)A_i=k_1(f_i,\cdot)\int_{D} k_2(x,y)A_i(y)dy,
\end{equation}
where $k_1$ and $k_2$ are scalar-valued kernel functions. By using neural networks $K_1,K_2$ to model $k_1,k_2$, we observe that $K_1$ and $K_2$ play the same role as the branch networks and the trunk networks in DeepONets. Next, we propose an architecture following an encoder-decoder framework. The neural networks $K_1,K_2$ encode the input to a $d$-dimensional vector field. We  proceed to equation \eqref{a4} in the $d$-dimensional vector field, and finally use a decoder to project back to the scaler field of interest.

We denote:
\begin{equation}
\begin{aligned}
a \odot b &= [a_1b_1,a_2b_2,\cdots,a_db_d]^T \\
a \otimes b &=[a_1b_1,a_1b_2,\cdots,a_2b_1,a_2b_2,\cdots,a_db_1,a_db_2,\cdots]^T,
\end{aligned}
\end{equation}
where $a=[a_1,a_2,\cdots,a_d]^T \in \mathbb{R}^d$ and $b=[b_1,b_2,\cdots,b_d]^T \in \mathbb{R}^d$.

Given data $\{f_i|_{P_k},u_i|_{P_k}\}|_{i=1}^n$ and $f$ in PDE \eqref{a1}, $u$ is the corresponding solution of $f$. Then we approximate $u$ by
\begin{equation}\label{a5}
v(x)=\sum_{i=1}^n K_1(f|_{P_k},f_i|_{P_k})\odot \int_D K_2(x,y) \otimes A_i(x)dy
\end{equation}
\begin{equation}
u^{app}(x)=Wv(x)+b,
\end{equation}
where $W \in \mathbb{R}^{1 \times d}, b \in \mathbb{R}$, and $K_1, K_2$ project the input to $\mathbb{R}^d$. Fig. \ref{fig:d} shows the architecture of our model.

It is worth noting that in practice, the integral in \eqref{a5} can not be calculated analytically, and we often represent it by numerical integration. We use the Gauss-Chebyshev integral, and the calculation   becomes
\begin{equation}
v(x)=\sum_{i=1}^n K_1(f|_{P_k},f_i|_{P_k}) \odot \sum_{y \in S_G} \omega_y K_2(x,y) \otimes A_i(x)
\end{equation}
\begin{equation}
u^{app}(x)=Wv(x)+b,
\end{equation}
where $S_G \subset D$ consists of fixed integration points determined by the Gauss-Chebyshev integral method, and $\omega_y$ are corresponding coefficients.

\section{Numerical Examples}\label{sec4}
We evaluate our model on several typical PDE datasets. As mentioned above, we solve the PDEs in  high resolution by a traditional numerical method and obtain our training data by downsampling. The numerical results show that our method can efficiently solve PDEs based on learning the operator. We demonstrate   our model's good generalization properties even with a small amount of training data, and show that the proposed method exhibits mesh-independent properties,   is trained at low resolutions, and can generalize at high resolutions.

Unless stated otherwise, in our model, the architectures of neural networks $K_1$ and $K_2$ are 128-256-256-128 and 128-128-128-128, respectively, with a ReLU activation function, and they encode inputs to a 32-dimensional vector field. $K_1$ is a convolutional neural network, and we first use one layer to enhance the input channel to 32 without an activation function. The kernel size is taken as 32, and we use padding to keep the channel unchanged. $K_2$ is a feedforward neural network (FNN). We choose $m=10$ for Chebyshev polynomials, $|S_G|=10$ in one-dimensional PDEs, and $|S_G|=100$ in two-dimensional PDEs for Gauss-Chebyshev integration. For all examples, we train our model in  4000 epochs using the Adam optimizer with an initial learning rate (LR) of 0.0002. Applying  polynomial decay to the LR, we set the power of the polynomial to $0.5$, and the LR will  decay to $2 \times 10^{-8}$ in 200 epochs. To efficiently train our model, we  increase the LR after each 200 epochs, and then decrease it. Due to the structural similarity, we choose \textbf{DeepONets} as our benchmark model. Meanwhile, we also take the traditional neural networks into account. \textbf{NN}: a point-wise feed-forward neural network. \textbf{CNN}: a fully convolutional neural network.

The PDEs we consider include the one-dimensional advection equation, burgers’ equation,   KdV equation, and   two-dimensional poisson equation. Advection   is one of the most important processes in atmospheric motion, and advection terms are included in the motion equation, heat flux equation, and water vapor equation in the atmospheric motion equations. Burgers' equation is a fundamental PDE in various fields of applied mathematics, involving many fields of applied mathematics such as fluid mechanics, nonlinear acoustics, and gas dynamics. The KdV equation, as  proposed by Dutch mathematicians Korteweg and de Vries,  describes nonlinear shallow water waves in fluid mechanics. The Poisson equation is commonly found in electrostatics, mechanical engineering, and theoretical physics.

\subsection{Advection Equation}\label{sec4.1}
We first consider the one-dimensional advection equation,
\begin{equation}
\begin{aligned}
\frac{\partial u}{\partial t}+\frac{\partial u}{\partial x} &=0,x \in (0,1),t \in (0,1] \\
u(0,t) &= u(1,t).
\end{aligned}
\end{equation}
Our task is to use the initial conditions $u(x,0)$ to give the solution at $u(x,t)$ at time t=1. We use Gaussian random fields to generate initial conditions with periodic boundary conditions,   numerically solve the equation at high resolution, and downsample the final training data   from the high-resolution numerical solution. The Gaussian random field is chosen as $\mathcal{N}(0,625(-\triangle+5^2I)^{-4})$, and we use Chebfun to solve the equation in $2^{11}$ resolution.
\begin{figure}[htb]
\centering
\subfloat[]{\includegraphics[width=.5\textwidth]{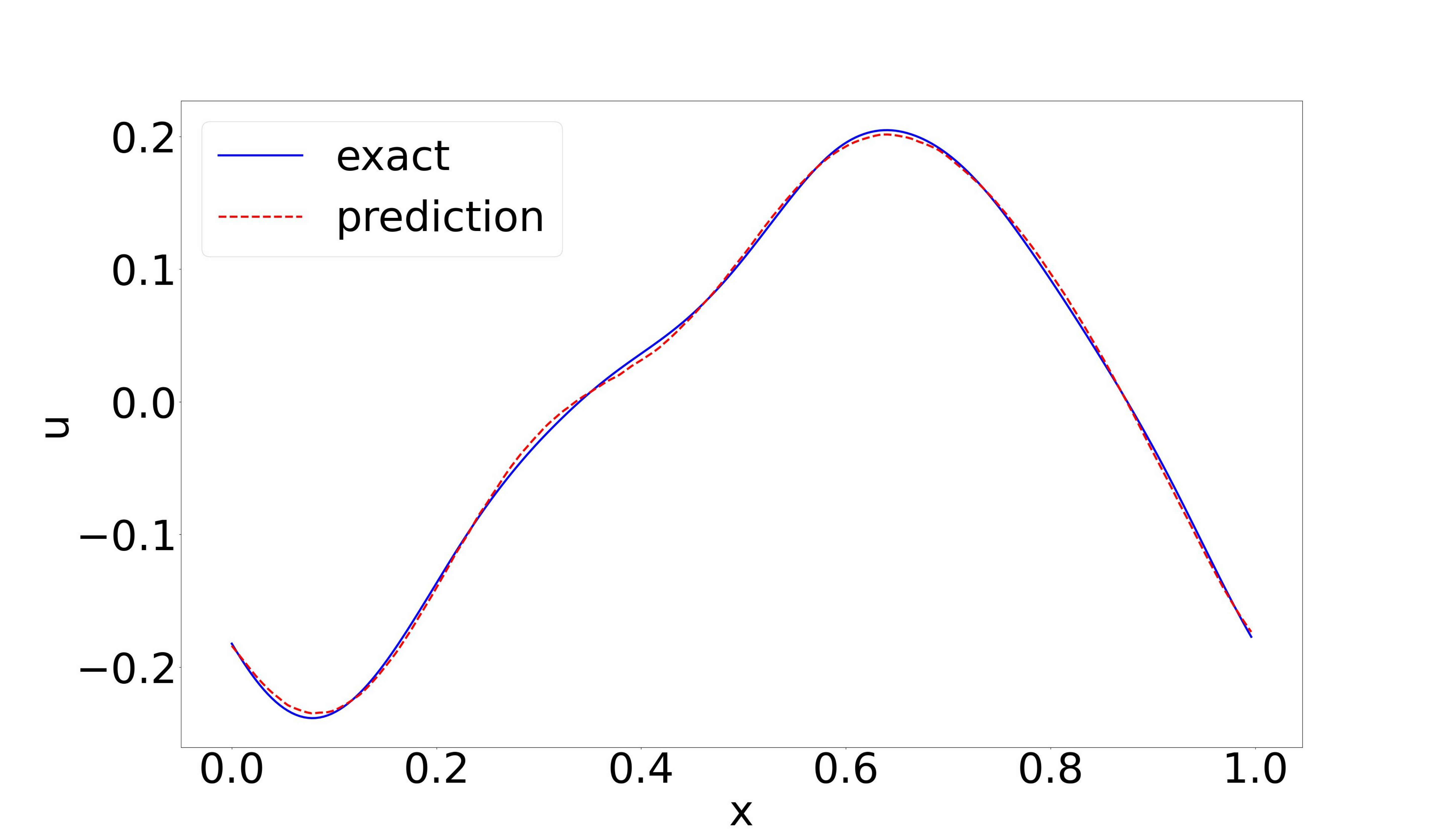}}
\subfloat[]{\includegraphics[width=.5\textwidth]{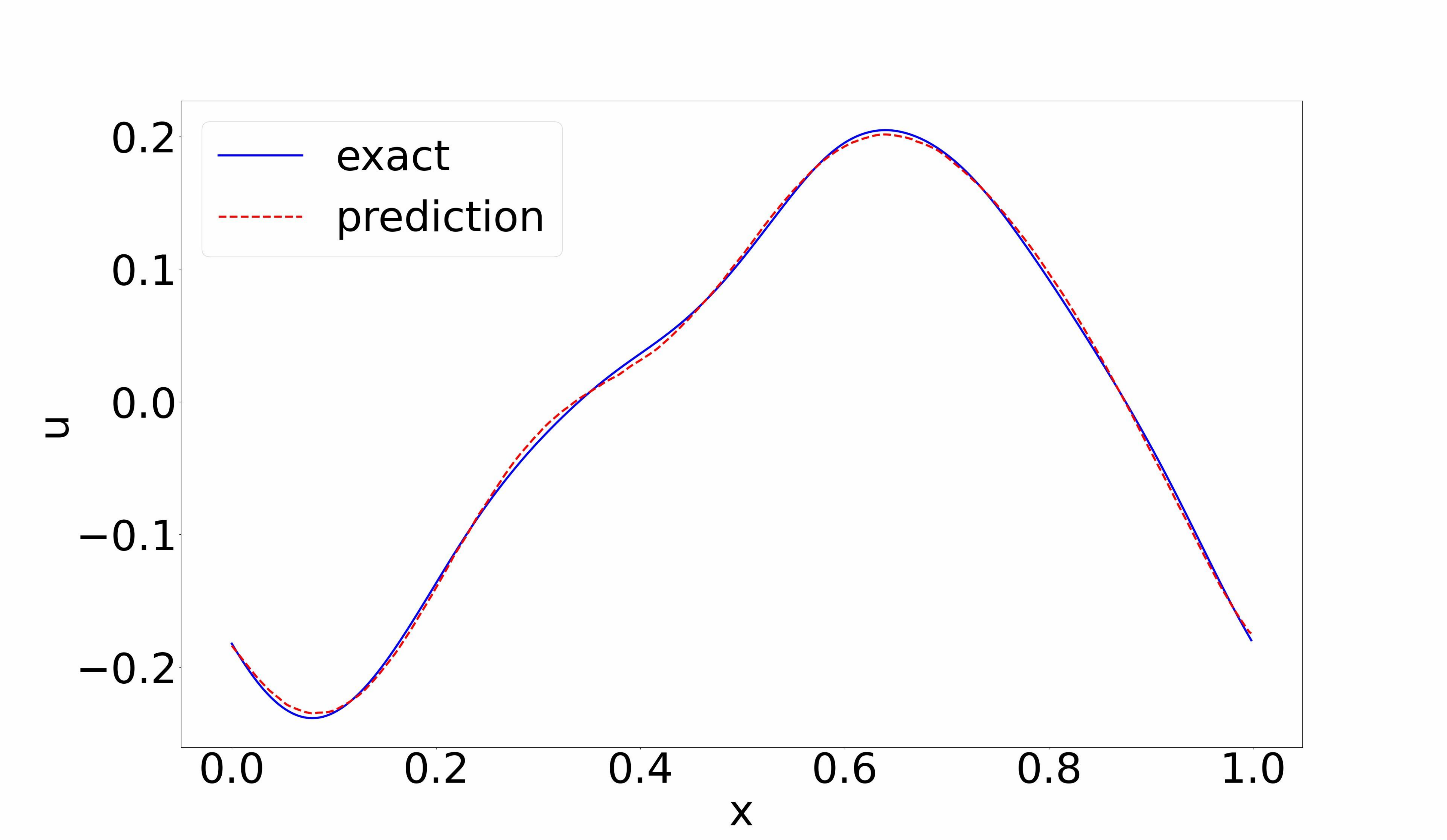}}
\caption{Comparison of  true solution and  learned model solution of advection equation at different resolutions: (a) $2^8$; (b)  $2^9$.}
\label{fig:1d}
\end{figure}

\begin{figure}[htb]
\centering
\includegraphics[width=.5\textwidth]{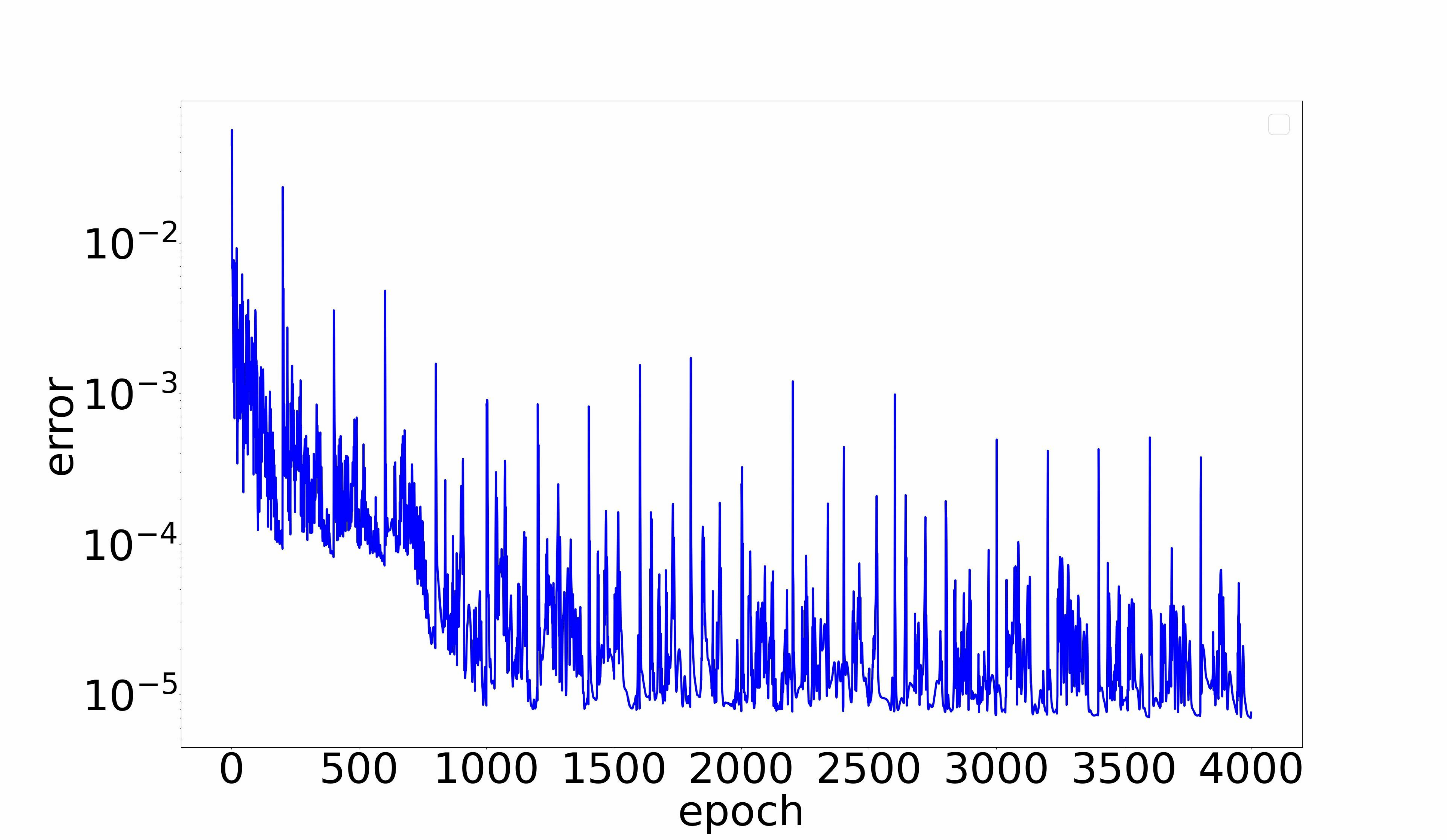}
\caption{Training loss history of advection equation.}
\label{fig:2d}
\end{figure}

The experimental results of the advection equation are shown in Fig. \ref{fig:1d} and Fig. \ref{fig:2d}. Given a new initial condition,  Fig. \ref{fig:1d}   compares the approximate solution generated by our model with the exact solution in a different resolution. Fig. \ref{fig:2d} shows the training error history in the view of $L_2$ error.

Our model is trained in  $2^8$ resolution and the result in Fig. \ref{fig:1d} (a) shows that the predicted solution fits accurately the exact solution. In addition, we evaluate our model at $2^9$ resolution. Fig. \ref{fig:1d} (b)  validates that the proposed method is able to find the solution of a high-resolution input after learning from lower-resolution data, which correctly learns the operator about the PDE. Due to the scheme, which we choose for LR, the variation trend of loss has large fluctuations. Meanwhile, it is seen from Fig. \ref{fig:2d} that after 1000 epochs, the training loss is able to reach $10^{-5}$ periodically. Consequently, the trained model demonstrates good generalization capability. Moreover, Our approach can achieve \textbf{1.7\%} relative $L_2$ error on the test dataset.

\subsection{Burgers’ Equation}\label{sec4.2}
The one-dimensional Burgers’ equation is a  PDE whose applications include modeling the one-dimensional flow of a viscous fluid. It takes the form
\begin{equation}
\begin{aligned}
\frac{\partial u}{\partial t}+u\frac{\partial u}{\partial x} &= v\frac{\partial^2 u}{\partial x^2},x \in (0,1),t \in (0,1]\\
u_0(x) &= u(x,t=0).
\end{aligned}
\end{equation}
The task for the neural operator is to learn the mapping of initial condition $u(x,t=0)$ to the solutions at $u(x,t=1)$. The initial condition $u_0(x)$ is generated according to $u_0 \sim \mu$, where $\mu=\mathcal{N}(0,625(-\triangle+25I)^{-4})$ with periodic boundary conditions. We set the viscosity to $v=0.01$ and solve the equation using a split-step method, where the heat equation part is solved exactly in Fourier space, and the  part is advanced again in Fourier space, using a very fine forward Euler method. We solve on a spatial mesh with resolution $2^{12}$, and use this dataset to downsample other resolutions.
\begin{figure}[htb]
\centering
\subfloat[]{\includegraphics[width=.5\textwidth]{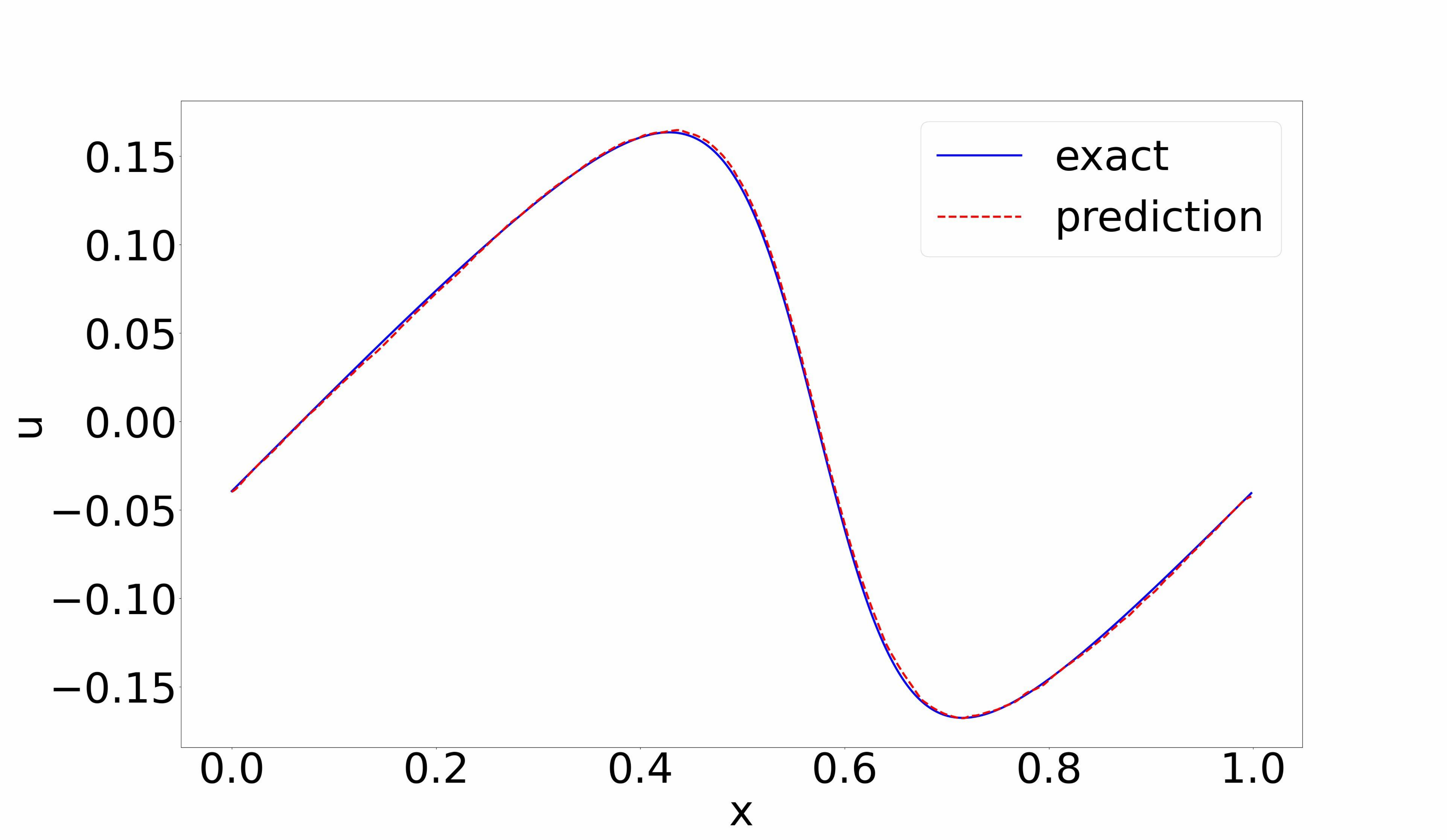}}
\subfloat[]{\includegraphics[width=.5\textwidth]{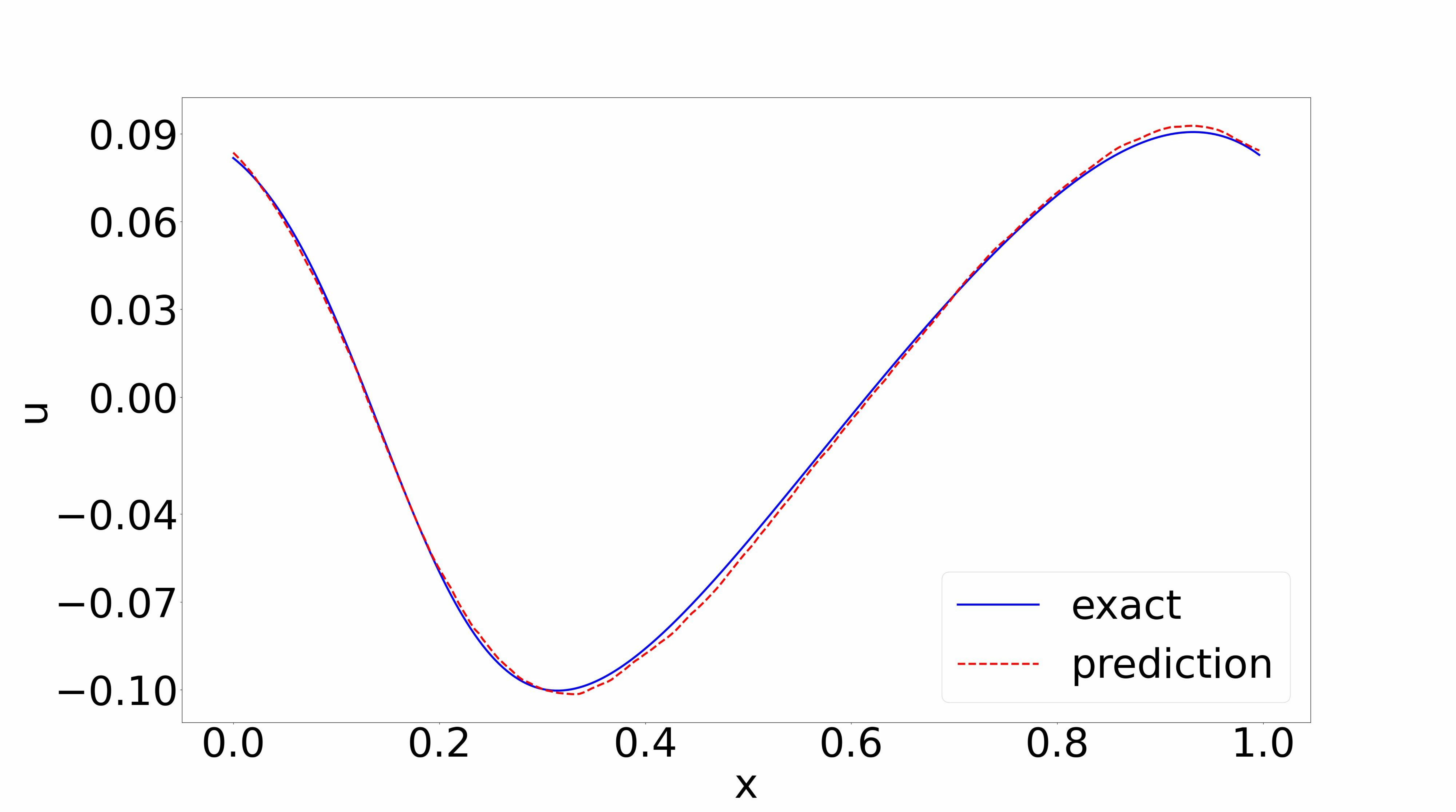}}
\caption{Burgers’ equation: Comparison of true solution and learned model solution in two inputs.}
\label{fig:3d}
\end{figure}

\begin{figure}[htb]
\centering
\includegraphics[width=.5\textwidth]{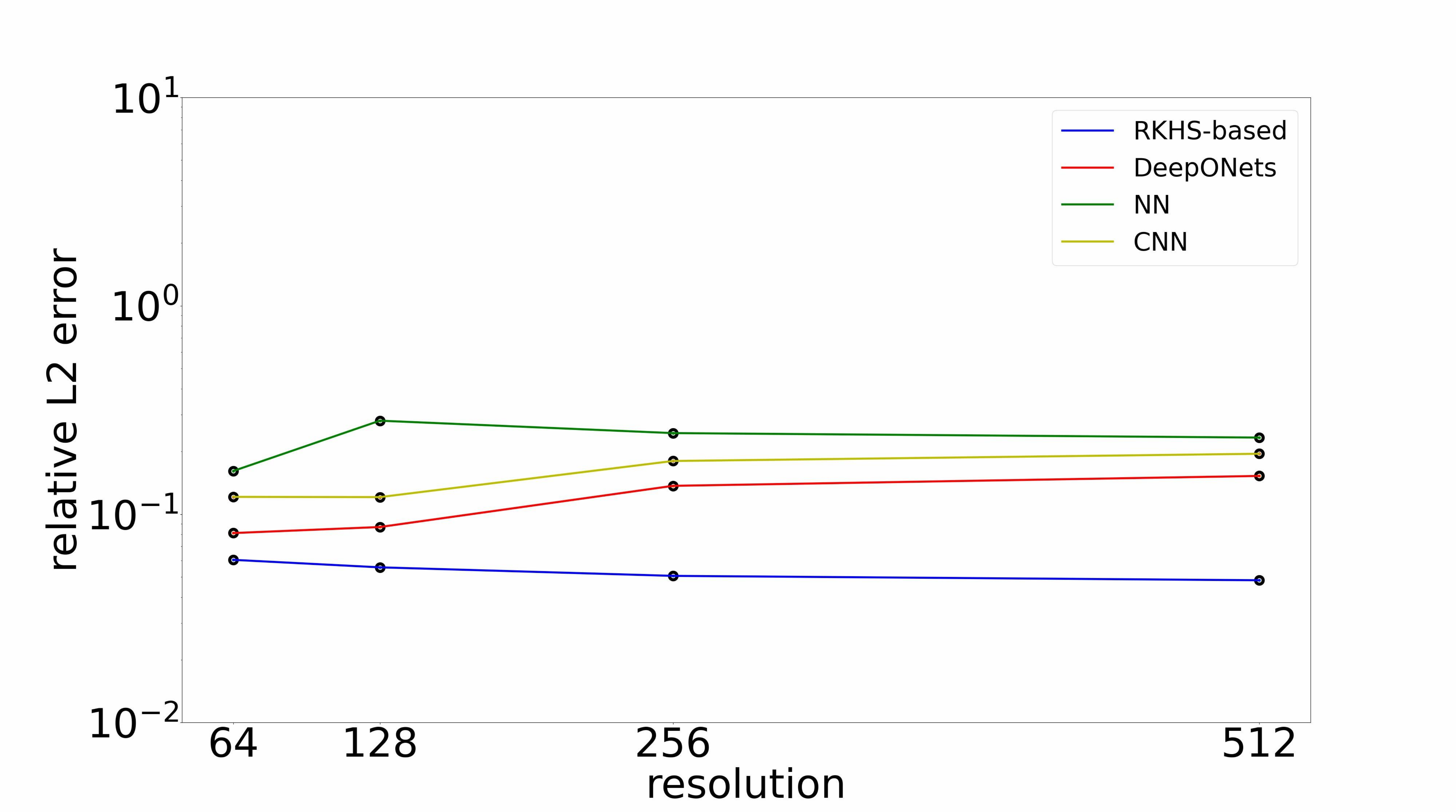}
\caption{Burgers’ equation validation at various input resolutions.}
\label{fig:4d}
\end{figure}

The results of our experiments are shown in Fig. \ref{fig:3d}, Fig. \ref{fig:4d}, Table \ref{tab:1}, and Table \ref{tab:2}. In Fig. \ref{fig:3d}, given new initial conditions, we compare the approximate solution generated by our model with the exact solution. Fig. \ref{fig:3d} and Table \ref{tab:1} compare our model and benchmark models in view of the relative $L_2$ error. Table \ref{tab:2} shows the results of our model trained in some resolutions and generalized at various input resolution  s.

The results show that the predicted solution by our approach simulates exact solution accurately, and obtains the lowest relative error compared to any of the benchmarks. Further, with the increasing resolution, we see that the relative error of our model gradually decreases, which learns better in  high resolution than in low resolution, while other models increase with resolution. Our model also demonstrates the property of mesh-independence that when the model is trained in a special resolution, the relative error is invariant with the testing resolution.
\begin{table}[htb]
\centering
\begin{tabular}{ccccc}
\toprule
\multicolumn{1}{c}{Resolutions} &\multicolumn{1}{c}{$s=64$} &\multicolumn{1}{c}{$s=128$} &\multicolumn{1}{c}{$s=256$} &\multicolumn{1}{c}{$s=512$} \\
\hline
\multicolumn{1}{c}{NN}& 0.1612  & 0.2810  & 0.2450 & 0.2332\\
\multicolumn{1}{c}{CNN}& 0.1211  & 0.1208  & 0.1800 & 0.1950\\
\multicolumn{1}{c}{DeepONets}& 0.0812  & 0.0868  & 0.1368 & 0.1526\\
\multicolumn{1}{c}{RKHS-based}& \textbf{0.0604}  & \textbf{0.0556}  & \textbf{0.0506}  & \textbf{0.0482}\\
\bottomrule
\end{tabular}
\caption{Benchmarks on Burgers’ equation at various input resolutions.}
\label{tab:1}
\end{table}

\begin{table}[htb]
\centering
\begin{tabular}{ccccc}
\toprule
\multicolumn{1}{c}{Resolutions} &\multicolumn{1}{c}{$s'=64$} &\multicolumn{1}{c}{$s'=128$} &\multicolumn{1}{c}{$s'=256$} &\multicolumn{1}{c}{$s'=512$} \\
\hline
\multicolumn{1}{c}{$s=64$}& 0.0604  & 0.0608  & 0.0610  & 0.0611\\
\multicolumn{1}{c}{$s=128$}& 0.0557  & 0.0556  & 0.0556  & 0.0556\\
\multicolumn{1}{c}{$s=256$}& 0.0506  & 0.0505  & 0.0506  & 0.0506\\
\multicolumn{1}{c}{$s=512$}& \textbf{0.0484}  & \textbf{0.0483}  & \textbf{0.0484}  & \textbf{0.0482}\\	
\bottomrule
\end{tabular}
\caption{Burgers’ equation: resolutions in training and testing.}
\label{tab:2}
\end{table}

\subsection{Korteweg-de Vries (KdV) Equation}\label{sec4.3}
The one-dimensional KdV equation takes the form
\begin{equation}
\begin{aligned}
\frac{\partial u}{\partial t}&=-0.5u\frac{\partial u}{\partial x}-\frac{\partial^3 u}{\partial x^3} x \in (0,1),t \in (0,1]\\
u_0(x) &= u(x,t=0).
\end{aligned}
\end{equation}
As mentioned above, we still solve the equation by learning the mapping of initial condition $u(x,t=0)$ to the solutions $u(x,t)$ at $t=1$. The initial condition $u_0(x)$ is generated according to a random field \cite{19} with a fluctuating parameter $\lambda$, and we take $\lambda=0.25$. We solve on a spatial mesh with resolution $2^{12}$, and use this dataset to downsample other resolutions to train our model.
\begin{figure}[htb]
\centering
\subfloat[]{\includegraphics[width=.5\textwidth]{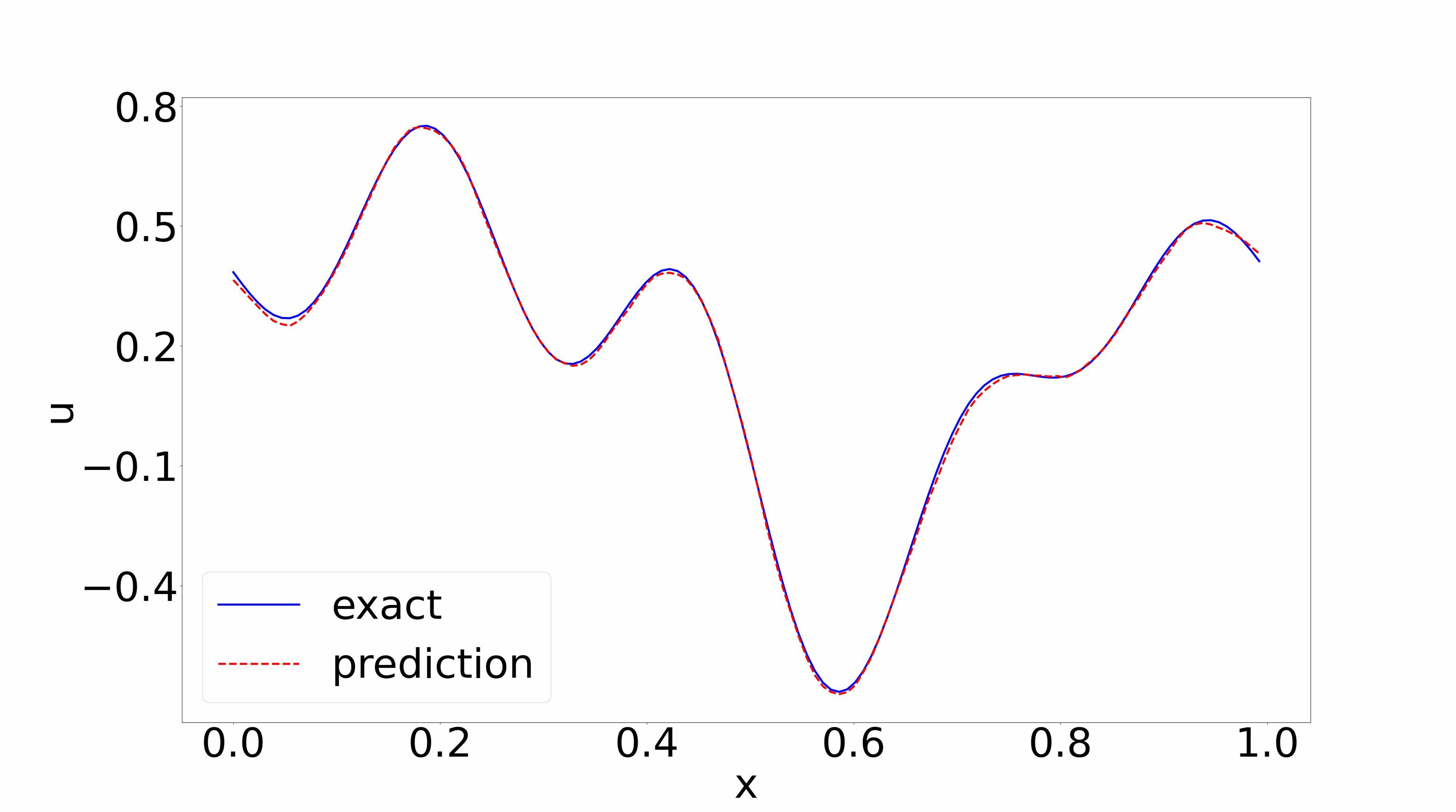}}
\subfloat[]{\includegraphics[width=.5\textwidth]{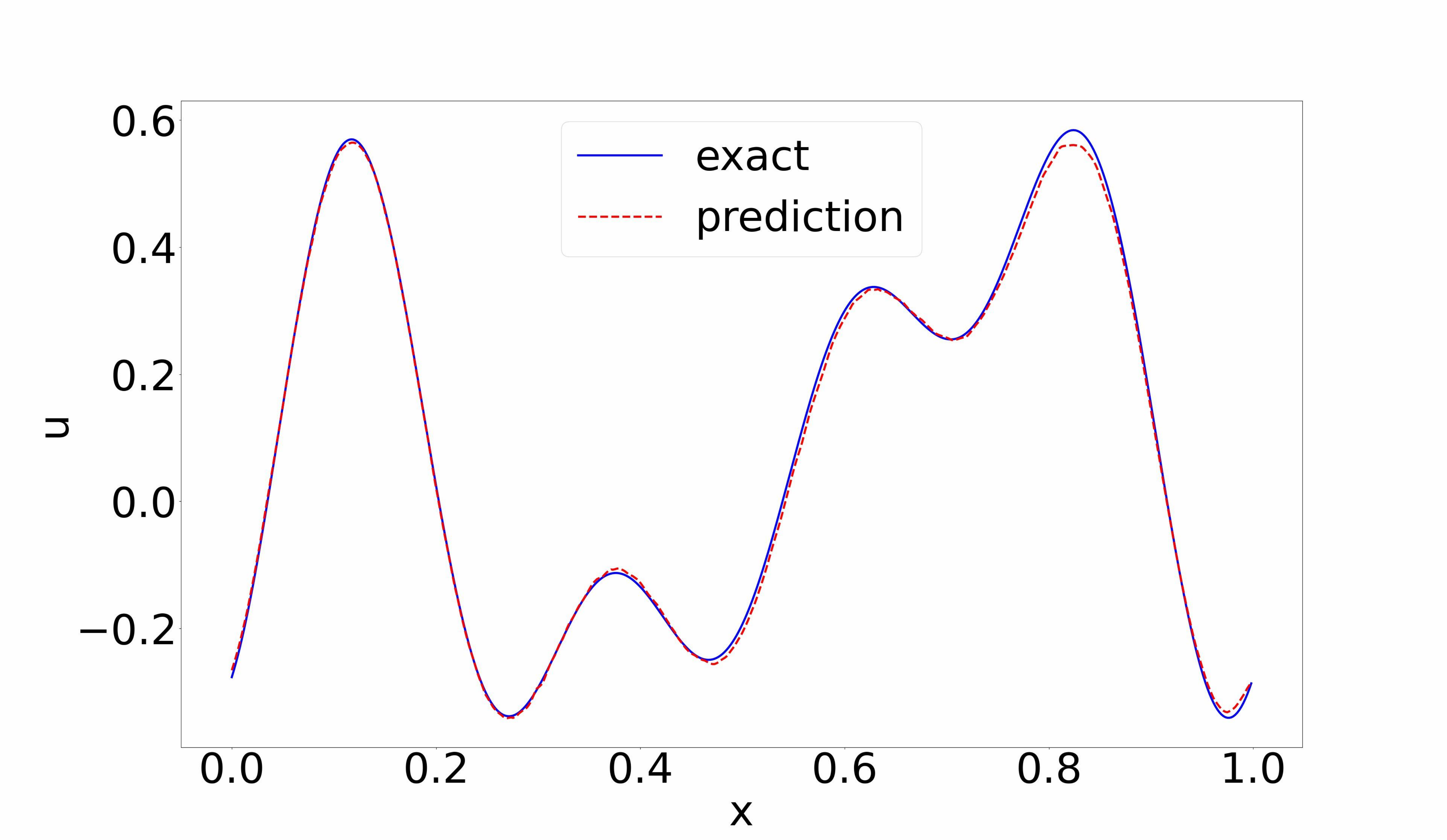}}
\caption{KdV equation: comparison of true solution and   learned model solution in two different inputs.}
\label{fig:5d}
\end{figure}

\begin{figure}[htb]
\centering
\includegraphics[width=.5\textwidth]{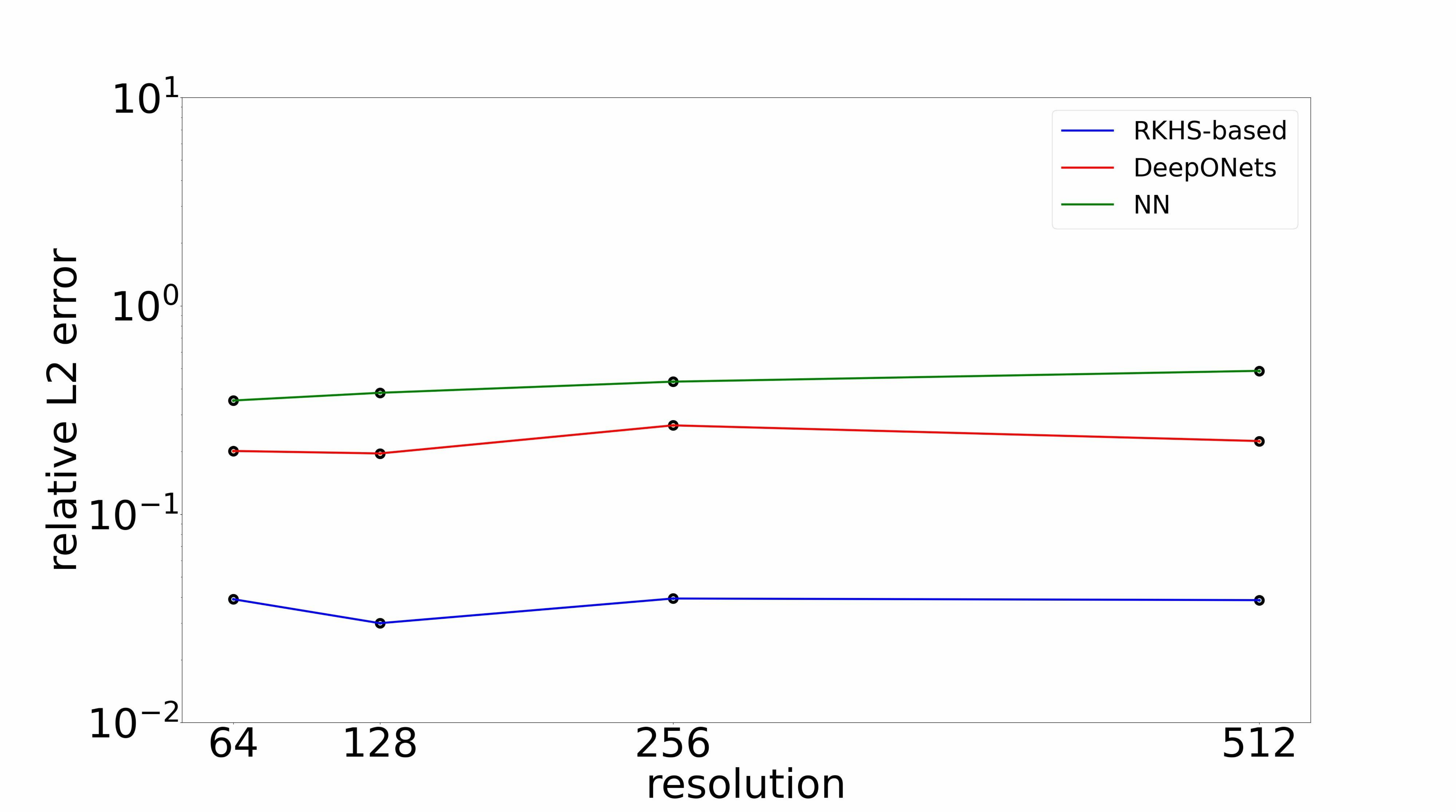}
\caption{KdV equation validation at various input resolutions.}
\label{fig:6d}
\end{figure}

Fig. \ref{fig:5d}, Fig. \ref{fig:6d}, Table \ref{tab:4}, and Table \ref{tab:3} show the results of our experiments. In Fig. \ref{fig:5d}, we compare the approximate solution generated by our model with the exact solution, given a new initial condition. Fig. \ref{fig:6d} and Table \ref{tab:4}   compare the relative $L_2$ error for our model and benchmark models. Table \ref{tab:3} shows the results for our model trained in some resolutions and generalized at various input resolutions.

It is worth noting that our model also obtains the lowest relative error compared to any of the benchmarks, and still maintains a small relative error compared to Burgers‘ equation. The predicted solution  fits accurately the exact solution, while the benchmarks  perform poorly for large fluctuations and even the CNN does not converge at all. At various input resolutions, our approach is able to achieve the same degree of relative error. Meanwhile, we again observe the invariance of the error with respect to the resolution. These results show that our approach is  mesh-independent  and has the property of stability.
\begin{table}[htb]
\centering
\begin{tabular}{ccccc}
\toprule
\multicolumn{1}{c}{Resolutions} &\multicolumn{1}{c}{$s=64$} &\multicolumn{1}{c}{$s=128$} &\multicolumn{1}{c}{$s=256$} &\multicolumn{1}{c}{$s=512$} \\
\hline
\multicolumn{1}{c}{NN}& 0.3512  & 0.3824  & 0.4325 & 0.4874\\
\multicolumn{1}{c}{CNN}& -  & -  & - & -\\
\multicolumn{1}{c}{DeepONets}& 0.2011  & 0.1957  & 0.2669 & 0.2244\\
\multicolumn{1}{c}{RKHS-based}& \textbf{0.0391}  & \textbf{0.0300}  & \textbf{0.0394}  & \textbf{0.0387}\\
\bottomrule
\end{tabular}
\caption{Benchmarks on KdV equation at various input resolution.}
\label{tab:4}
\end{table}

\begin{table}[htb]
\centering
\begin{tabular}{ccccc}
\toprule
\multicolumn{1}{c}{Resolutions} &\multicolumn{1}{c}{$s'=64$} &\multicolumn{1}{c}{$s'=128$} &\multicolumn{1}{c}{$s'=256$} &\multicolumn{1}{c}{$s'=512$} \\
\hline
\multicolumn{1}{c}{$s=64$}& 0.0391  & 0.0399  & 0.0404  & 0.0406\\
\multicolumn{1}{c}{$s=128$}& \textbf{0.0291}  & \textbf{0.0300}  & \textbf{0.0310}  & \textbf{0.0315}\\
\multicolumn{1}{c}{$s=256$}& 0.0396  & 0.0395  & 0.0394  & 0.0394\\
\multicolumn{1}{c}{$s=512$}& 0.0386  & 0.0388  & 0.0387  & 0.0387\\	
\bottomrule
\end{tabular}
\caption{KdV equation: Resolutions in training and testing.}
\label{tab:3}
\end{table}

\subsection{Poisson equation}\label{sec4.4}
We consider the poisson equation in the case with a small amount of data \cite{20},
\begin{equation}
\begin{aligned}
-\bigtriangleup u(x)&= f(x),x \in D \\
u(x)&=0,x\in \partial D.
\end{aligned}
\end{equation}

Consider the 2D case, in which the source function $f$ is generated by $-a(x^2-x+y^2-y)$, i.e.,
\begin{equation}
\begin{aligned}
\partial_{xx}u+\partial_{yy}u&=-a(x^2-x+y^2-y),(x,y) \in D\\
u&=0,(x,y)\in \partial D,
\end{aligned}
\end{equation}
where $D=[0,1]^2$, and the constant $a$ controls the source term. Obviously, the analytical solution is $u=\frac{a}{2}x(x-1)y(y-1)$.

\begin{figure}[htb]
\centering
\subfloat[]{\includegraphics[width=.5\textwidth]{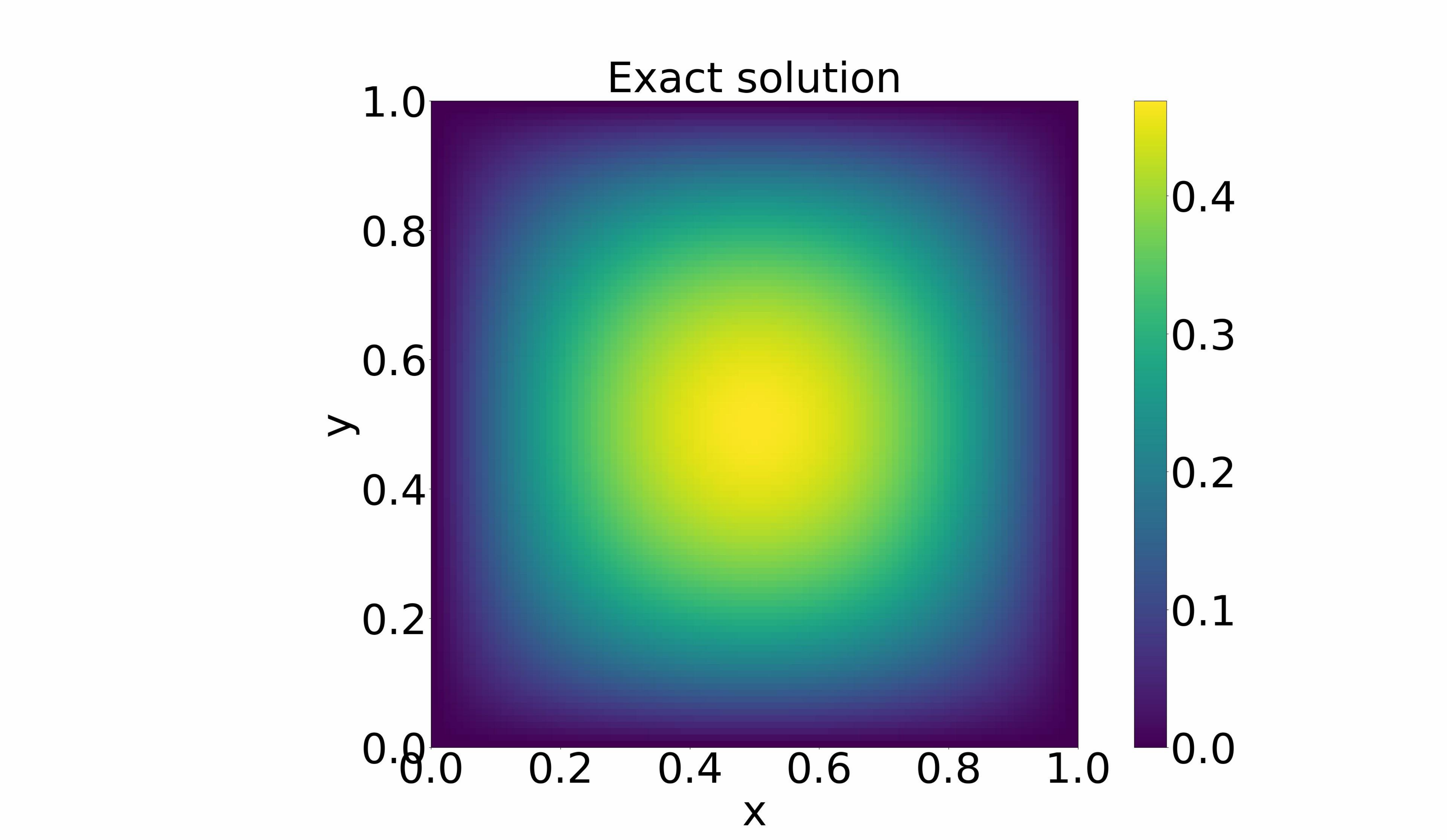}}
\subfloat[]{\includegraphics[width=.5\textwidth]{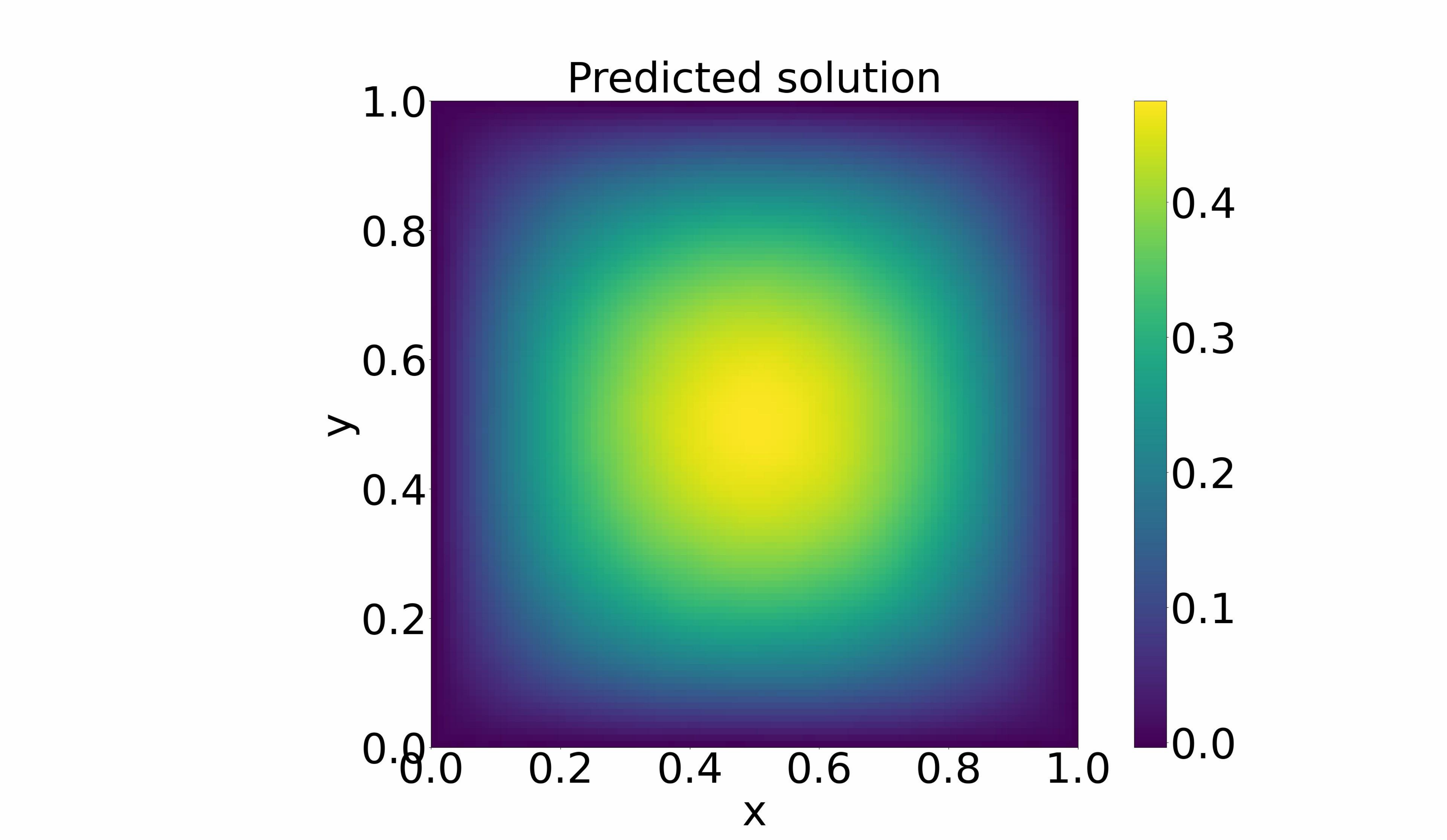}}
\\
\centering
\subfloat[]{\includegraphics[width=.5\textwidth]{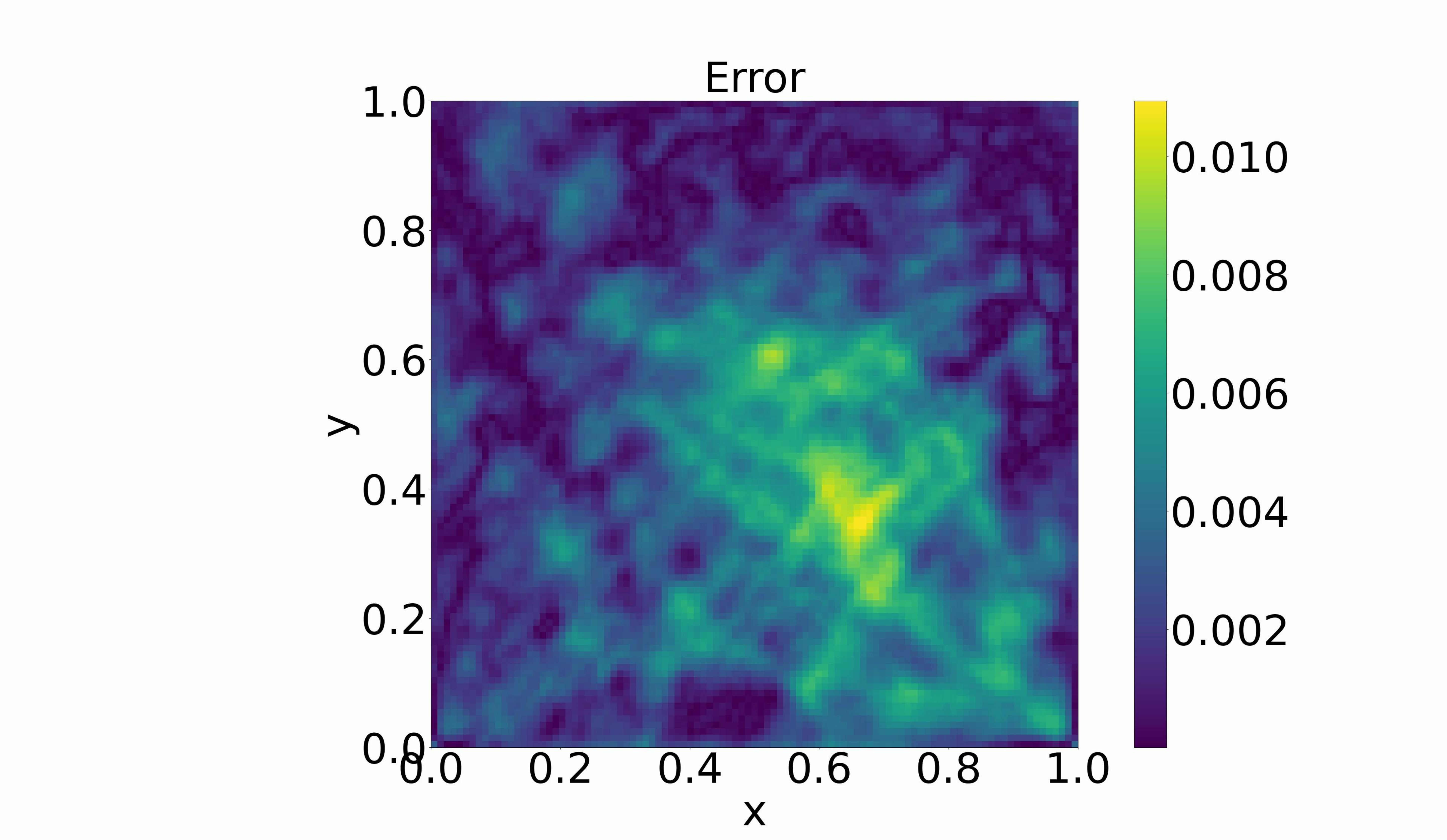}}
\caption{Poisson equation: comparison between exact   and predicted solutions on 101×101 grid points corresponding to source terms determined by $a=15$: (a) exact solution; (b) predicted solution; (c) difference between exact and predicted solutions. }
\label{fig:7d}
\end{figure}

In this experiment, the settings of the training process are different from those mentioned above. During training, we choose a set of source functions, i.e., we sample $10$ source functions $f$ from a selected region, which is implemented by sampling control parameter $a$ uniformly from $\{10k\}_{k=1}^{20}$, and take $P_k$ as 64 × 64 isometric grid points. For each epoch, we randomly sample 200 points in $[0,1]^2$ and sample 100 points in each line of the boundary. It is worth noting that in previous experiments, the training points were also $P_k$. Finally, we test the performance of our well-trained model on $a=15$. To visualize its performance, we show the exact  and   predicted solutions on 101×101 isometric grid points. To more intuitively compare these   solutions, we calculate their difference   at each point and indicate the error by colors, as shown in Fig. \ref{fig:7d}. The relative $L_2$ error is \textbf{1.5\%}.

We also test DeepONets on a two-dimensional Poisson equation with the similar training process, as we sample $K=10$ source functions $f$ from the selected region and evaluate them on 64 × 64 isometric grid points. For each epoch, we randomly sample 200 points in $[0,1]^2$ and 100 points in each line of the boundary. We observe that the testing relative $L_2$ error is \textbf{62.3\%}. Moreover, we evaluate them on 101 × 101 isometric grid points and also use the grid points as training points to train Deeponets, NN and CNN. The testing relative $L_2$ error is as follows, \textbf{Deeponets}: \textbf{0.26\%}, \textbf{NN}: $4.6\%$ and \textbf{CNN}: $8.6\%$. The results show that when the data is sufficient, Deeponets achieves a better performance than other models. However, in the case of a small amount of data, our model may be a better choice.

\section{Conclusion}\label{sec5}
In this study, we transformed the solution of PDEs to the problem of data-driven learning of the operator that maps between two function spaces. Motivated by the fundamental properties of the function-valued RKHS, we proposed an architecture following the form given by the representer theorem, whose  kernel was constructed by a Hilbert-Schmidt integral operator. We showed that the proposed model   performed well on several kinds of PDEs, even  solutions with  high fluctuation. Numerical experiments showed that our approach was  mesh-independent, with the ability to find the solution of a high-resolution input after learning from lower-resolution data, and had the property of stability. Another advantage of our model was that it only used a small amount of data to train, but had  good generalization.


\bibliographystyle{elsarticle-num}
\bibliography{bkj}

\end{document}